\documentstyle{article}
\newcommand{\C}{{{\rm C}\!\!\!{\rm I}\,}}
\newcommand{\R}{{{\rm I}\!{\rm R}}}
\newcommand{\Z}{{{\rm Z}\!\!{\rm Z}}}
\newcommand{\eqd}{\buildrel {\rm def}\over =}

\newcommand{\im}{{\rm {Im\ }}}
\newcommand{\sign}{{\rm {sign\ }}}
\renewcommand{\L}{{\cal L}}
\newcommand{\A}{{\cal A}}
\renewcommand{\hat}{\widehat}
\renewcommand{\tilde}{\widetilde}
\newtheorem{theorem}{Theorem}
\begin{document}
\begin{large}

\title{Entire periodic functions with plane zeros}
\author{L.I.Ronkin
\and A.M.Russakovskii}
\date{}
{\maketitle}

\begin{abstract}
We give a complete description of divisors of entire periodic
functions in $\C^n$ with plane zeros.
\end{abstract}

\section {Introduction}

The present note is close by the considered problems to the
papers \cite {Ronkin-divisors,Ronkin-algebraic} and continues
investigations of entire periodic functions of several complex
variables and their divisors started in these papers.

In \cite{Ronkin-divisors} conditions were studied, under which
an $n$ - periodic positive divisor in $\C^n$ is a divisor of
some entire $n$ - periodic function. In \cite {Ronkin-algebraic}
a description was given, of those polynomials in $n$ variables
whose divisors after their $n$ - periodic "reproduction"
generate $n$ - periodic divisors in $\C^n.$ Here we give a
complete description of divisors of entire periodic functions
with plane zeros\footnote {The term "function with plane zeros"
is utilized for functions whose zero set is a union of
hyperplanes. Such functions have been object of study in
a number of papers, particularly in 
\cite {Gruman,Sekerin,Papush,Papush-Russakovskii}}.  Besides that we give a
more thorough analysis of the conditions on divisors introduced
in \cite{Ronkin-divisors}.

We denote further by $Z$ an $n$ - periodic divisor with
linearly independent periods
$\omega_1\in\R^n,\ldots,\omega_n\in\R^n$ and by $f(z)$ an
arbitrary entire function in $\C^n$ whose divisor $Z_f$
coincides with $Z.$ In this situation,
$$f(z+\omega_p) = e^{g_{\omega_p(z)}}f(z),\quad p=1,\ldots, n,$$
where the entire functions $g_p=g_{\omega_p}$ are defined by
$f$ up to constants of the form $2\pi im, \ m\in \Z.$

Denote by $\Delta_p=\Delta_{\omega_p}$ the difference operator
$$\Delta_p h =h(z+\omega_p) - h(z).$$
The integer quantities
$$N_{pq}=N(\omega_p, \omega_q : Z)\eqd\frac{1}{2\pi
i}(\Delta_p g_q - \Delta_q g_p)$$
do not depend on $z$ nor on the choice of $f(z)$ (see next
section). The skew-symmetric matrix
$$\tilde{N}_Z=\{N_{pq}\}_{p,q=1}^n$$
will be called the {\it index} of the periodic\footnote{For
brevity, in what follows we will write simply
"periodic divisor" instead of "$n$ - periodic divisor with
periods $\omega_1,\ldots,\omega_n\in\R^n$".} divisor $Z.$

Without loss of generality we assume henceforth that
$$\omega_1=(1,0,\ldots,0),\ldots, \omega_n=(0,\ldots,0,1).$$

Using the concept of the index of a periodic divisor, one can
reformulate theorem 1 from \cite{Ronkin-divisors} in the
following way.

\medskip

{\bf Theorem A.}
\begin{it}
A periodic divisor $Z\subset\C^n$ is a divisor
of some entire periodic function if and only if
$$\tilde{N}_Z=0.$$
\end{it}

\medskip

In \cite {Ronkin-divisors}, a geometric sufficient condition
for a divisor to be a divisor of a periodic entire function was
also given:

\medskip

{\bf Theorem B.}
\begin{it}
Let $\alpha_{pq}$ be the automorphism of $\C^n$ permuting the
$z_p$ and $z_q$ coordinates.

If a periodic divisor $Z\subset \C^n$ is symmetric, i.e.
invariant with respect to $\alpha_{pq},\ \forall p,q,$ then
there exists an entire periodic function $F(z)$ with $Z_F=Z.$
\end{it}

\bigskip

In the next section, we make some elementary analysis of the
properties of an index of a periodic divisor which allows to
give a somewhat more general sufficient condition.

\bigskip

\begin{theorem}
Let $\beta_j$ be the automorphism of $\C^n$ taking $z_j$ to
$-z_j$ and not changing the rest of coordinates and let
$\alpha_{pq}$ be the same as in theorem B.

If for some splitting of the set $\{1,\ldots,n\}$ into
nonintersecting parts $I$ and $J$ a divisor $Z\subset \C^n$ is
invariant with respect to $\beta_k, \ \forall k\in I,$ and with
respect to $\alpha_{pq},\ \forall p,q\in J,$ then there exists
an entire periodic function $F(z)$ with $Z_F=Z.$
\end{theorem}

\bigskip

The main subject of our study is periodic divisors with plane
components, i.e. periodic divisors in $\C^n$ whose carrier is a
union of hyperplanes. To formulate the main result of the
paper, introduce additional denotions.

Let
$$l(z)=l_{a,c} (z)=<a,z>+c,\quad a\in\C^n,\ c\in\C$$
and let
$$L=L_{a,c}=\{z\in\C^n:\ l(z)=0\}$$
be the corresponding hyperplane. Denote by $\hat{S}_L$ the
periodic "reproduction" of the hyperplane $L,$ i.e.
$$\hat{S}_L=\bigcup\limits_{k\in\Z^n}\{z\in\C^n:\quad
l(z-k)=0\}.$$

Note that this operation does not necessarily produce an
analytic set. If, however, it {\it is} analytic, we identify it
with the divisor having $\hat{S}_L$ as its carrier with
multiplicity 1 on it.

Consider the set $$A_l=\{x\in\R^n:\ <a,x>=0\}.$$ This is a
subspace in $\R^n.$ Denote by $m$ the dimension of its
orthogonal complement in $\R^n.$ Obviously, either $m=1$ or
$m=2.$ Denote the set of those linear functions $l(z)$ for
which $m=1$ by ${\cal L}_1$ and the set of those ones for which
$m=2$ by ${\cal L}_2.$

If $l_{a,c}$ is given, consider a "parallelogram"
$$P_{pq}=\{w\in\C:\ w=\alpha a_p + \beta a_q; \ 0\leq \alpha<1,
\ 0\leq \beta <1\}.$$
If
$$a_q \neq 0\ {\rm {and}}\ \im \frac{a_p}{a_q}=0,$$
then $P_{pq}$ has no interior and becomes an interval.

Denote by $\nu_{pq}$ the number of points $w\in\C$ of the form
$w=<k,a>,\ k\in\Z^n$ belonging to $P_{pq}.$ If $a_p$ or $a_q$
equals $0,$ we set $\nu_{pq}=0.$

If $l\in\L_2$ and $\hat{S}_L$ is a divisor, then, as we will
see in section 3, $$\im \frac{a_p}{a_q}\neq 0,$$ (and hence
$P_{pq}$ is nondegenerate) for at least one pair $(p,q).$

\bigskip

\begin{theorem}
A periodic divisor $Z\subset \C^n$ is a divisor of some entire
periodic function with plane zeros if and only if it may be
represented in the form
$$Z=Z'+Z'',$$
where
$$Z'=\sum\limits_{j=1}^{\mu_1} \hat{S}_{L_j},\quad \mu_1
\leq \infty, \quad l_j=<a^{(j)},z>+c_j,\quad l_j\in{\cal L}_1$$
and
$$\lim_{j\to\infty}\frac{|\im
c_j|}{|a^{(j)}|}=\infty,\quad \rm{if}\ \mu_1=\infty,$$
and where
$$Z''=\sum\limits_{j=1}^{\mu_2} \hat{S}_{L_j},\quad \mu_2
< \infty, \quad l_j=<a^{(j)},z>+c_j, \quad l_j\in{\cal
L}_2$$
and
$$\sum\limits_{j=1}^{\mu_2}\nu^{(j)}_{pq} \sign \im
\frac{a_q^{(j)}}{a_p^{(j)}},\quad \forall p,q \ {\rm{such\
that}}\ a_p^{(j)}\neq 0.$$
\end{theorem}

\bigskip

{\bf Remark.} In section 3 where this theorem is proved we give
also the explicit form of the corresponding entire function.

\bigskip

\section {Properties of the index and proof of theorem 1}

In \cite{Ronkin-divisors} it was noted that the quantity
$N_{pq}$ does not depend on $z$ and is an integer. It is easy to
see also that this quantity does not depend on the choice of
the function $f$ defining the divisor $Z.$ Indeed, if $f$ and
$\tilde{f}$ are entire functions with
$$Z_f=Z_{\tilde{f}}=Z,$$
then
$$f=\tilde{f}e^h$$
with $h$ entire, and hence
$$g_p=\tilde{g}_p +\Delta_p h.$$
The function $\tilde{g}_p$ here is defined similar to $g,$ i.e.
by the equality
$$\tilde{f}(z+\omega_p)=\tilde{f}(z)e^{\tilde{g}_p(z)}.$$

We have then
$$\Delta_p g_q - \Delta_p \tilde{g}_q = \Delta_p \Delta_q h,$$
and due to commutativity of the operators $\Delta_p$ and
$\Delta_q$ it follows that
$$\Delta_p g_q - \Delta_q g_p = \Delta_p \tilde{g}_q - \Delta_q
\tilde{g}_p.$$

We outline below several elementary properties of the matrix
$\tilde{N}_Z.$

\bigskip

1) if $k\in\Z,\ k>0,$ then
$$N(k\omega_p, \omega_q: Z) = kN(\omega_p, \omega_q: Z).$$

According to the definition, we have to determine functions
$g_q(z)=g_{\omega_q}(z)$ and $g_p(z)=g_{k\omega_p}(z)$ first.
If $f(z)$ is an entire function with $Z_f=Z,$ then $g_p(z)$ is
defined by the relation
$$f(z+k\omega_p)=f(z)e^{g_p(z)}.$$
Hence one easily sees that we may take
$$g_{k\omega_p} (z) = \sum\limits_{j=0}^{k-1} g_{\omega_p}
(z+j\omega_p).$$
Similarly, the operator $\Delta_p=\Delta_{k\omega_p}$
applied to a function $h$ gives
$$[\Delta_{k\omega_p} h]\ (z)=\sum\limits_{j=0}^{k-1}
[\Delta_{\omega_p} h]\ (z+j\omega_p).$$

Therefore we have
$$2\pi i N(k\omega_p, \omega_q:
Z) =  \Delta_{k\omega_p}g_q - \Delta_{\omega_q} g_p$$
$$= \sum\limits_{j=0}^{k-1}
\left[ [\Delta_{\omega_p}g_q](z+j\omega_p) -
[\Delta_{\omega_p} g_{\omega_p}] (z+ j\omega_p)\right]$$
$$=  k\cdot 2\pi iN(\omega_p, \omega_q: Z).$$

\bigskip

2) If a divisor $Z^*$ is obtained from the divisor $Z$ with the
help of the mapping
$$\beta_p:z\mapsto (z_1,\ldots, z_{p-1}, -z_p,
z_{p+1},\ldots,z_n),$$
then
$$N(\omega_p, \omega_q: Z^*) = - N(\omega_p, \omega_q: Z).$$

To show that this is true, note first that if $Z=Z_f$ for some
entire function $f$ then $Z^*=Z(f^*),$ where $f^*(z)=f(\beta_p
z).$ Therefore the corresponding functions $g_j^*$ for $j\neq p$
are related to $g_j$ by the equalities $g_j^*=g_j(\beta_p z).$
For $j=p$ one has
$$e^{g^*_p(z)} = \frac{f^*(z+\omega_p)}{f^*(z)}
= \frac{f(\beta_p z+ \beta_p \omega_p)}{f(\beta_p z)}
= \frac{f(\beta_p z-\omega_p)}{f(\beta_p z)}
= e^{-g_p (\beta_p z -\omega_p)}$$
and hence one can take
$$g^*_p(z)=-g_p (\beta_p z -\omega_p).$$

>From the above relations between $g_j^*$ and $g_j$ it follows
that

$$2\pi i N(\omega_p, \omega_q:Z^*) = \Delta_p g^*_q -
\Delta_q g_p^*$$
$$= [g_q (\beta_p z+\beta_p\omega_p)-g_q(\beta_p z)] -
[-g_p(\beta_p z -\omega_p+\omega_q) + g(\beta_p z -\omega_p)]$$
$$ = -[g_q (\beta_p z) - g_q (\beta_p z-\omega_p)] + [g_p
(\beta_p z -\omega_p+\omega_q) - g(\beta_p z -\omega_p)]$$
$$= (\Delta_q g_p - \Delta_p g_q)|_{\beta_p z - \omega_p}
= -2\pi i N(\omega_p, \omega_q: Z).$$

\bigskip

3) If a divisor $Z^*$ is obtained from the divisor $Z$ with the
help of the mapping $\alpha_{pq}:\C^n\to\C^n$ permuting the
$z_p$ and $z_q$ coordinates and leaving the other coordinates
unchanged, then
$$N(\omega_p, \omega_q:Z^*) = - N(\omega_p, \omega_q: Z).$$

To prove this, consider as before functions $f, g_p$ and $f^*,
g^*_p$ corresponding to the divisors $Z$ and $Z^*.$ Evidently,
one can take $f^*(z)=f(\alpha_{pq}(z).$ Then
$$\frac{f^*(z+\omega_p)}{f^*(z)}=\frac{f(\alpha_{pq}z+
\alpha_{pq}\omega_p)}{f(\alpha_{pq}z)}
=\frac{f(\alpha_{pq}z+\omega_q)}{f(\alpha_{pq}z)},$$
whence one can take
$$g^*_p=g_q(\alpha_{pq}z).$$

This implies

$$2\pi i N(\omega_p, \omega_q: Z^*) =
[g_q^*(z+\omega_p)-g_q^*(z)] - [g^*_p(z+\omega_q)-g^*_p (z)]$$
$$= [g_p(\alpha_{pq}z +\omega_q) - g_p(\alpha_{pq}z)]-
[g_q(\alpha_{pq}z +\omega_p)-g_q(\alpha_{pq}z)]$$
$$= (\Delta_q g_p - \Delta_p g_q)|_{\alpha_{pq}z}
= 2\pi i N(\omega_q, \omega_p: Z)
= - 2\pi i N(\omega_p, \omega_q: Z)$$

\bigskip

4) If a divisor $Z$ is a sum of periodic divisors $Z_1$ and
$Z_2$ then
$$\tilde{N}_Z=\tilde{N}_{Z_1} + \tilde{N}_{Z_2}.$$

This property is obvious.

\bigskip

We would like to note that the theorem B mentioned in the
introduction follows immediately from theorem A and the
property 3). In \cite{Ronkin-divisors} it was proved (also
 deduced from theorem A) in a less elementary way.

If one uses the property 2) together with the property 3), then
theorem A yields immediately our theorem 1.

\section {Linear periodic divisors. Proof of theorem 2}

Next we consider periodic divisors with plane components, i.e.
divisors whose carrier is the union of hyperplanes of the form
$\{z:\ <a,z>+c=0\}.$

Let, as before,
$$l(z)=l_{a,c} (z)=<a,z>+c,\quad a\in\C^n,\ c\in\C$$
and let
$$A=A_l=\{x\in\R^n:\ l(z-x)\equiv l(z)\}=\{x\in\R^n:\
<a,x>=0\},$$
$$A^\perp=A_l^\perp = \{\xi\in \R^n:\ <\xi, x>=0,\ \forall x\in
A\}.$$
Denote the dimension of $A^\perp$ by $m.$ Let $\Lambda_1,
\ldots, \Lambda_m$ be a basis in $A^\perp$ and
$\Lambda_{m+1},\ldots, \Lambda_n$ be a basis in $A.$ Then the
vectors $\Lambda_1,\ldots,\Lambda_n$ form a basis in both
$\R^n$ and $\C^n$ (over $\R$ and $\C$ respectively).
Hence one can represent $z\in \C^n$ as follows:
$$z=\sum\limits_{j=1}^n \zeta_j(z) \Lambda_j,$$
where $\zeta_j(z)=<\Omega_j, z>$ and $\Omega_j$ is the $j$ - th
row of the matrix $\Omega$ reciprocal to the matrix
$\Lambda$ with columns $\Lambda_1,\ldots,\Lambda_n.$

>From the definition of the spaces $A^\perp$ and $A$ it follows
that
$$l(z)=c+<a,\zeta_1(z)\Lambda_1+\ldots+\zeta_m(z)\Lambda_m>=
l^*(\zeta')|_{\zeta'=\zeta'(z)},$$
where
$$\zeta'=(\zeta_1,\ldots,\zeta_m),\quad
l^*(\zeta')=\sum\limits_{j=1}^m b_j\zeta_j,\quad b_j =
<a,\Lambda_j>.$$

Let $L$ be the hyperplane $\{z:\ l(z)=0\}.$ Consider
the periodic "reproduction" of the hyperplane $L,$
i.e.  $$\hat{S}_L=\bigcup\limits_{k\in\Z^n}\{z\in\C^n:\quad
l(z-k)=0\}.$$

This set is obviously periodic and if it is analytic,
it can be considered as a periodic divisor, i.e. the carrier of
the divisor is $\hat{S}_L$ and the multiplicity on it equals 1.

>From the results in \cite{Ronkin-algebraic} it follows in
particular that $\hat{S}_L$ will be a divisor if and only if
the space $A^\perp$ possesses a basis consisting of vectors
with integer coordinates and the function $l^*(\zeta')$ is
hypoelliptic, that is, satisfies the condition:
{\it
the distance
from a point $\xi'\in\R^m$ to the set $\{\zeta'\in\C^m:\
l^*(\zeta')=0\}$ tends to $\infty$ when $\xi'$ tends to
$\infty.$}

Let $\hat{S}_L$ be a divisor. We are interested in the question
when this divisor is a divisor of some entire periodic
function.

We have already split all linear functions $l(z)$
into two sets $\L_1$ and $\L_2$ with regard to the
dimension $m$ of the space $A^\perp.$

Consider the case $l\in\L_1$ first. In this case
$$A^\perp=\{x\in\R^n:\ x=ta^0,\ t\in \R\},\quad a^0\in\R^n,\
a=w_0a^0,$$
where $w_0\in\C,$ and if $\hat{S}_L$ is a divisor, without loss
of generality we may assume that $a=a^0=k^0\in\Z^n,$ and that
the coordinates of $k^0$ have no common divisor.

Put
$$f_l(z)=\sin[\pi(<k^0,z>+c)]e^{i\pi (<k^0, z> +c)}.$$

Evidently, $f_l$ is an entire periodic function with divisor
$$Z_f=\bigcup\limits_{j=-\infty}^{\infty}\{z:\ <k^0,z>+c=j\} =
\bigcup\limits_{k\in\Z^n} \{z:\ <k^0, z>+c=<k^0,k>\} =
\hat{S}_L.$$

Hence we get an affirmative answer to our question when
$l\in\L_1.$

\bigskip

Let now $l\in\L_2.$ Then
$$l^*(\zeta')=b_1\zeta_1 + b_2\zeta_2+c,$$
and
$$l(z)=b_1\zeta_1(z) + b_2\zeta_2(z)+c,$$
where
$$b_1=<a,\Lambda_1>,\quad b_2=<a,\Lambda_2>,\quad
\zeta_1(z)=<\Omega_1, z>,\quad \zeta_2(z)=<\Omega_2, z>.$$

If $\hat{S}_L$ is a divisor, then the function $l^*(\zeta')$ is
hypoelliptic, and we may assume that the elements of the matrix
$\Lambda$ are integers. It is immediate to check that the
function $l^*(\zeta_1, \zeta_2)$ is hypoelliptic if and only if
$$b_1\neq 0,\ b_2\neq 0,\ {\rm {and}}\ \im
\frac{b_2}{b_1}\neq 0,$$ or equivalently, when $\displaystyle
\im\frac{a_p}{a_q}\neq 0$ for at least one pair of indices $p$
and $q.$

\bigskip

First we consider the simplest case with regard to the number
of variables, $n.$ Note that if $n=2,$ the original function
$l(z_1, z_2)$ is also hypoelliptic.  To compute the index of
the divisor $Z=\hat{S}_L$ in this case, we construct an
entire function of one complex variable $w$ which has simple
zeros precisely at points $w=p+qT,\ (p,q)\in\Z^2$ where
$T=\frac{a_2}{a_1}.$

Suppose that $\im T>0$ and consider the product
$$\Phi_T(w)=\prod\limits_{q=0}^\infty\{2ie^{-i\pi(w-qT)}\sin\pi(w-qT)\}
\prod\limits_{q=-\infty}^{-1}\{2ie^{i\pi(w-qT)}\sin\pi(w-qT)\}\eqno (1).$$

Since
$$2ie^{-iw}\sin w=1+O(e^{-2|\im w|})\quad {\rm {when}\ } \im
w\to-\infty,$$
and
$$2ie^{iw}\sin w=1+O(e^{-2|\im w|})\quad {\rm {when}\ } \im
w\to+\infty,$$
the above product converges and hence is an entire function.

Similarly, if $\im T<0$ we set
$$\Phi_T(w)
= \prod\limits_{q=0}^\infty\{2ie^{i\pi(w-qT)}\sin\pi(w-qT)\}
\prod\limits_{q=-\infty}^{-1}\{2ie^{-i\pi(w-qT)}\sin\pi(w-qT)\}\eqno (2).$$

Next we denote
$$u^\pm_q(w)=2ie^{\pm i\pi(w-qT)}\sin\pi(w-qT).$$
It is obvious that
$$u^\pm_q (w+1)=u^\pm_q(w) \quad {\rm {and}\ }\quad
u^\pm_q(w+T)=u^\pm_{q-1}(w).$$
Hence
$$\Phi_T(w+1)=\Phi_T(w)\eqno (3)$$
and for $\im T>0$
$$\Phi_T(w+T)= \prod\limits_{q=0}^\infty u^-_q(w+T)
\prod\limits_{q=-\infty}^{-1}u^+_q(w+T)
= \frac{u_{-1}^-(w)}{u_{-1}^+(w)} \Phi_T(w) = e^{-2\pi i
(w+T)}\Phi_T(w).\eqno (4)$$
Similarly, when $\im T<0,$
$$\Phi_T(w+T)=e^{2\pi i (w+T)} \Phi_T(w).\eqno (5)$$

Consider the function
$$f(z_1, z_2)= f_{a,c}(z_1,
z_2)=\Phi_T(z_1+Tz_2+\frac{c}{a_1}).\eqno(6)$$

This entire function vanishes if and only if
$$z_1+Tz_2+\frac{c}{a_1}=p+qT,\quad (p,q)\in\Z^2,$$
that is, when
$$l(z_1-p, z_2-q)=0,\quad (p,q)\in\Z^2.$$
Thus $Z_f=\hat{S}_L.$

>From (3)-(6) it follows that the corresponding functions $g_1$
and $g_2$ may be taken in the following form:
$$g_1\equiv 0$$
and
$$g_2(z) = (-\sign \im T) \cdot 2\pi i
(z_1+Tz_2+\frac{c}{a_1}+T).$$

Therefore
$$N(\omega_1, \omega_2:\hat{S}_L) = - \sign \im T,$$
and hence the divisor $\hat{S}_L$ for $l\in\L_2$ cannot be a
divisor of any entire periodic function.

\bigskip

Assume now that the dimension $n>2,$ and let $l\in\L_2.$
Suppose $\im \frac{a_p}{a_q}\neq 0$ for some $p$ and $q.$

Fix these $p$ and $q.$ To compute $N_{pq},$ we will construct
a function $f(z)$ having divisor $Z_f=\hat{S}_L$ in the form
$$f(z)=\Psi(\frac{l(z)}{a_p}),$$
where the entire function $\Psi(w),\ w\in\C,$ periodic with
period 1, is to be defined below.

Let $T_j=\frac{a_j}{a_p},\ j=1,\ldots, n.$
In order that $Z_f=\hat{S}_L,$ zeros of the function
$\Psi(w)$ must be simple and have the form
$$\{w\in\C:\ w=k_p+\sum\limits_{j=1}^n k_j T_j\}_{k=(k_1,
\ldots,k_n)\in\Z^n},$$
i.e. they must make up a periodic set $Z_\Psi$ with $n$ periods
$T_1,\ldots,T_n.$

Denote by $\A_l$ the set
$$\A_l=\{w\in\C^n:\ w=<a,k>,\ k\in\Z^n\}$$
and let $x_1=0,\ldots, x_{\nu_{pq}}$ be all the points from
$\A_l$ belonging to the parallelogram $$P_{pq}=\{w:\ w=\alpha
a_p + \beta a_q,\ 0\leq \alpha < 1, 0\leq \beta <1\}.$$
Then
$$\A_l=\bigcup\limits_{j=1}^{\nu_{pq}} A_j,$$ where
$$A_j=\{w\in\C:\ w=k_p a_p + k_q a_q + x_j, \ (k_p,
k_q)\in\Z^2\}.$$

>From all above it follows that $Z_\Psi$ may be represented in
the form of a union of $\nu_{pq}$ identical
nonintersecting sets
$$Z_j=\{w\in\C: a_p w\in A_j\}$$ each having periods 1 and
$T_q,$ and being shifted relative to each other:
$$Z_\Psi=\bigcup\limits_{j=1}^{\nu_{pq}}
(Z_1+\frac{x_j}{a_p}).$$

Using the functions $\Phi_T(w)$ constructed above, we set
$$\Psi (w)=\bigcup\limits_{j=1}^{\nu_{pq}}
\Phi_{T_q}(w-\frac{x_j}{a_p}).$$

By virtue of our construction the function
$$f(z)=\Psi(\frac{l(z)}{a_p})$$ has divisor $\hat{S}_L.$

Since the function $\Psi$ is periodic with period 1,
one can take
$$g_p(z)=\log \frac{f(z+\omega_p)}{f(z)}\equiv 0,$$
and
$$g_q(z)=-\sign \im T_q\cdot 2\pi i
\sum\limits_{j=1}^{\nu_{pq}}
(\frac{l(z)}{a_p}-\frac{x_j}{a_p}+T_q),\eqno (7)$$ which
implies that $$N(\omega_p,
\omega_q:\hat{S}_L)=-\nu_{pq}\cdot\sign\im T_q\neq 0.$$

Therefore the divisor $\hat{S}_L$ cannot be a divisor of
any entire periodic function.

\bigskip

Let, as above, $l\in\L_2$ and assume that $\hat{S}_L$ is a
divisor. In what follows we will need to compute $N_{pq}$ also
for such $p$ and $q$ that $\im \frac{a_p}{a_q} = 0.$ Note that
in this case $\frac{a_q}{a_p}$ is a rational number since in
the two-dimensional (over $\Z$) set $\{\frac{a_j}{a_p},\
j=1,\ldots,n\}$ there are already two $\Z$ - independent
periods ($1$ and $\frac{a_i}{a_p}\not\in \R$ for some $i$).
So, there are two integers, $k$ and $j$ such that
$$j a_p = k a_q.$$

To compute $N_{pq}$ we use properties 1) and 2) of the index.
In view of these properties
$$N(\omega_p, \omega_q: Z)= \frac{1}{k} N(\omega_p,
k \omega_q: Z),$$
so it is enough to compute $N(\omega_p, k \omega_q:Z).$
By the above definition of $f(z),$
$$e^{g_{k\omega_q}(z)}=\frac{f(z+k\omega_q)}{f(z)} =
\frac{\Psi(\frac{l(z+k\omega_q)}{a_p})}{\Psi(\frac{l(z)}{a_p})}
=\frac{\Psi (\frac{l(z)}{a_p}+\frac{k
a_q}{a_p})}{\Psi(\frac{l(z)}{a_p})}=\frac{\Psi
(\frac{l(z)}{a_p}+j)}{\Psi (\frac{l(z)}{a_p})}=1,$$
since the function $\Psi(z)$ has period $1$ and $j\in\Z.$
Therefore we can take
$g_{k\omega_q}\equiv 0$
and since $g_p = g_{\omega_p}\equiv 0,$ it follows that
$N(\omega_p, k_2\omega_q: Z)=0,$ and consequently $N_{pq}=0.$

\bigskip

Now we have all the necessary information to describe the
structure of divisors $Z$ of entire periodic functions with
plane zeros. Such a divisor is obviously periodic and has only
hyperplane irreducible components. From the above
considerations it follows that a periodic divisor with plane
components may be represented in the form
$$Z=Z'+Z'',$$
where
$$Z'=\sum\limits_{j=1}^{\mu_1} \hat{S}_{L_j},\quad
l_j\in\L_1,\quad \mu_1\leq\infty,$$
and
$$Z'=\sum\limits_{j=1}^{\mu_2} \hat{S}_{L_j},\quad
l_j\in\L_2,\quad \mu_2\leq\infty,$$
so that each $\hat{S}_{L_j}$ is a divisor itself.

It is immediate to check that each function $l\in\L_2$ vanishes
at some point of $\R^n=\R^n+i0.$ Hence the corresponding set
$\hat{S}_L$ has nonempty intersection with the qube
$$\{x\in\R^n:\ |x_j|\leq 1,\ j=1,\ldots,n\}.$$
Since a divisor can have only a finite number of irreducible
components intersecting a fixed compact, it follows that
$\mu_2<\infty$ and hence
$$\tilde{N}_{Z''}=\sum\limits_{j=1}^{\mu_2}\tilde{N}_{Z_j},$$
where $Z_j=\hat{S}_{L_j}.$

>From theorem A and property 4) of indices it follows now that
the divisor $Z''$ will be a divisor of an entire periodic
function if and only if
$$ \sum\limits_{j=1}^{\mu_2}
N(\omega_p,\omega_q:\hat{S}_{L_j})=0,\quad \forall p,q.$$

In terms of coefficients ($l_j=<a^{(j)},z>+c_j$) the last
condition may be rewritten as
$$\sum\limits_{j=1}^{\mu_2} \nu_{pq}^{(j)}\sign (\im
\frac{a_q^{(j)}}{a_p^{(j)}})=0,\quad \forall
p,q=1,\ldots,n\ {\rm {with}}\ a_p^{(j)}\neq 0.\eqno(8)$$

The corresponding entire periodic function, i.e. entire
periodic function $F_2$ with $Z_{F_2}=Z''$ is constructed as
follows.

First, we make additional investigation of the case
$$l\in\L_2,\ \hat{S}_L\ {\rm {is\ a\ divisor}}.$$
As it has been noted before, in this case
$$l(z)=b_1\zeta_1(z)+ b_2\zeta_2(z),$$ where
$$b_1\neq 0,\ b_2\neq 0,\ \im \frac{b_2}{b_1}\neq 0, \
\zeta_1=<k,\Omega_1>, \ \zeta_2 =<k,\Omega_2>,$$
and $\Omega_1, \Omega_2$ are vectors with rational coordinates.
The set
$${\cal A}_l\eqd \{w=<a,k>,\ k\in\Z^n\}$$ has the form
$${\cal A}_l=\{w:\ w=<k, b_1\Omega_1 + b_2\Omega_2>\}.$$

Since $\Omega_1$ and $\Omega_2$ have rational coordinates,
the set ${\cal A}_l,$ which is obviously an additive group, has
no density points in $\C$ and is therefore a free group with 2
generators. Hence there exist such $$w_1=w_1(l)\in \C\setminus
0,\ w_2=w_2(l)\in\C\setminus 0,\ \frac{w_2}{w_1}=
T(l)\in\C\setminus\R,$$ that
$${\cal A}_l=\{w: w= j_1 w_1 + j_2 w_2,\ (j_1, j_2)\in\Z^2\}.$$

Particularly,
$$a_p= m_{1p}w_1 + m_{2p} w_2,\quad m_{ip}=m_{ip}(l)\in\Z,\
i=1, 2.$$

Set
$$F_{(l)}(z)= \Phi_T (\frac{l(z)}{w_1}),$$
where the function $\Phi$ is defined by (1) or (2) depending on
the sign of $\im T{\rm(}=T(l){\rm )}.$

Then
$$Z_{F_{(l)}}=\bigcup\limits_{(j_1, j_2)\in\Z^2} \{z:\ l(z)=j_1
w_1+j_2w_2\} = \bigcup\limits_{k\in\Z^n} \{z: l(z) = <a,k>\} =
\hat{S}_L,$$
and by (3) and (4)
$$F_{(l)} (z+\omega_p) =e^{\tilde{g}_{p,l} (z) } F_{(l)}(z),$$
with
$$\tilde{g}_{p,l}(z)=-2\pi i\sign T(l) \left[
\frac{m_{1p}l(z)}{w_1}+\frac{m_{1p}(m_{1p}+1)}{2} T(l)\right].$$

Getting back to our divisor $\displaystyle
Z''=\bigcup\limits_{j=1}^{\mu_2} \hat{S}_{L_j},$ set
$$\tilde{F}=\prod\limits_{j=1}^{\mu_2}F_{(l_j)}(z).\eqno(9)$$

>From all above it follows that
$$Z_{\tilde{F}}=Z''$$ and that
$$\tilde{F}(z+\omega_p)=e^{G_p(z)}\tilde{F}(z)$$
with
$$G_p(z)=2\pi i
\sum\limits_{j=1}^{\mu_2}
\left(
\sign\im T(l_j)\cdot
\left[
\frac{m_{1p}(l_j) l_j(z)}{w_1^{(j)}}+
\frac{m_{1p}(l_j)(m_{1p}(l_j)+1)}{2} T(l_j)\right]
\right).$$

Therefore the functions $G_p(z)$
are linear in all variables, i.e.  have the form
$$\sum\limits_{j=1}^n \sigma_j^p z_j + \tau_p.$$

Since $N(\omega_p, \omega_q: Z'')=0$ and since  $\Delta_q
G_p=\sigma_q^p,$ it follows that $\sigma_q^p=\sigma_p^q,\
\forall p, q.$

Put
$$H(z)=\frac{1}{2} \sum\limits_{j=1}^{n}\sigma_j^jz_j^2
+\sum\limits_{j=1, j\neq p}^n \sigma_j^p z_j z_p +
\sum\limits_{j=1}^{n} (\tau_j-\frac{1}{2}\sigma_j^j) z_j.\eqno 
(10)$$

Obviously,
$$\Delta_p H = G_p.$$

Therefore the function
$$F_2(z)=\tilde{F}(z)e^{-H(z)}\eqno (11)$$
will be the required periodic function with divisor $Z''.$

\bigskip

It remains to consider the sum $Z'.$ This sum will be a divisor
if and only if each compact in $\C^n$ is intersected only by a
finite number of divisors $\hat{S}_{L_j}.$ Since $l_j\in\L_1,$
we may assume that each function $l_j$ has the form
$$l_j=<k^{(j)},z>+c_j,\quad k^{(j)}\in\Z^n,\ c_j\in\C.$$
Hence the projection of the divisor $\hat{S}_{L_j}$ onto $\R^n$
has nonempty intersection with each cube
$$\{x\in\R^n:\ t_j\leq x_j\leq t_j+1,\ j=1,\ldots,n\},$$
and the projection of $\hat{S}_{L_j}$ onto $i\R^n$ is the plane
$$\{iy:\ y\in\R^n,\ <k^{(j)}, y>+ \im c_j =0\}.$$

Therefore the sum
$\displaystyle \sum\limits_{j=1}^{\mu_1}\hat{S}_{L_j}$ is a
divisor if and only if either $\mu_1<\infty$ or $\mu_1=\infty$
and
$$\lim_{j\to\infty}\frac{|\im c_j|}{|k^{(j)}|}=\infty.$$

It was shown above that each divisor $\hat{S}_{L_j}$ is a
divisor of some entire periodic function. So if $\mu_1<\infty,$
the divisor $Z'$ will also be a divisor of a periodic function
equal to the product of the corresponding functions for
$\hat{S}_{L_j}.$ We will show that a similar fact holds also
in the case $\mu_1=\infty.$

In order to prove this set
$$\gamma_j=\frac{\im c_j}{|k^{(j)}|}$$ and form the product
$$F_1(z)=\prod\limits_{j=1}^\infty \phi_j(z).\eqno (12)$$
with
$$\phi_j = 2 i\sin \pi (<k^{(j)},z>+c_j)\cdot e^{\sign \gamma_j
i\pi (<k^{(j)}, z>+c_j)}\eqno (13).$$

Taking into account the above mentioned asymptotics of
$2ie^w\sin w$ we obtain that for $|\im z|<R,$
$$F_1(z)=\prod\limits_{j=1}^\infty (1+O(e^{-|k^{(j)}|
(|\gamma_j|-R)}).\eqno (14)$$

By our assumption, $|\gamma_j|\to\infty$. Hence (14) implies
that the product $F_1$ converges.

Summarizing all above we come to the following conclusion:

\bigskip

{\bf Theorem 2.}
\begin{it}
A periodic divisor $Z\subset \C^n$ is a divisor of some entire
periodic function $F(z)$ with plane zeros if and only if it may
be represented in the form $$Z=Z'+Z'',$$ where
$$Z'=\sum\limits_{j=1}^{\mu_1} \hat{S}_{L_j},\quad \mu_1
\leq \infty, \quad l_j=<a^{(j)},z>+c_j,\quad l_j\in{\cal L}_1$$
and
$$\lim_{j\to\infty}\frac{|\im
c_j|}{|a^{(j)}|}=\infty,\quad \rm{if}\ \mu_1=\infty,$$
and where
$$Z''=\sum\limits_{j=1}^{\mu_2} \hat{S}_{L_j},\quad \mu_2
< \infty, \quad l_j=<a^{(j)},z>+c_j, \quad l_j\in{\cal
L}_2$$
and
$$\sum\limits_{j=1}^{\mu_2}\nu^{(j)}_{pq} \sign \im
\frac{a_q^{(j)}}{a_p^{(j)}},\quad \forall p,q \ {\rm{such\
that}}\ a_p^{(j)}\neq 0.$$

Under these conditions the function $F(z)$ is representable
in the form
$$F(z)=e^{\Psi(z)}F_1(z)F_2(z),$$
where $\Psi(z)$ is an arbitrary entire periodic function, and
the functions $F_1(z)$ and $F_2(z)$ are constructed with
regard to divisors $Z'$ and $Z''$ by the above mentioned
method, i.e. by virtue of (9) - (14). \end{it}

\obeylines
Theory of Functions Department
Mathematical Division
Institute for Low Temperature Physics
47 Lenin Avenue
310164 Kharkov
Ukraine

\end{large}

\bigskip

\begin{thebibliography}{MMMM}

\bibitem[Gr]{Gruman} L.Gruman. The regularity of growth
of entire functions whose zeros are hyperplanes. Ark. Mat. 10,
1972, p. 23-31.

\bibitem[Pa]{Papush}
D. E. Papush. Entire functions of
several variables with a regular set of 'plane' zeros. Sib.
Math. J., 32, 1991, p. 120-130.

\bibitem[PaRu]{Papush-Russakovskii} D.E.Papush,
A.M.Russakovskii. Interpolation on plane sets in $\C^2.$ Ann.
Fac. Sci. Toulouse, 1, n 3, 1992, p. 337 - 362.

\bibitem[Ro1]{Ronkin-divisors} L.I.Ronkin. Entire periodic
functions of several variables, periodic divisors. To appear in
Math. Physics, Math. Anal. and Geom.

\bibitem[Ro2]{Ronkin-algebraic} L.I.Ronkin. Periodic algebraic
divisors in $\C^n$. To appear in Math. Physics, Math. Anal. and
Geom.

\bibitem[Se]{Sekerin} A.B.Sekerin. On construction of entire
functions with given growth. Sib. Math. J., 27, n 3, 1986, p.
179 - 192.

\end{thebibliography}
\end{document}